\documentclass[10pt, twoside]{article}
\usepackage[b5j,hmargin={1 in,0.6 in},vmargin={1 in,0.6 in}]{geometry}
\usepackage{times}
\usepackage[labelsep=space, labelfont=bf]{caption}
\usepackage{amsmath,amssymb,amsthm,latexsym}
\makeatletter	
\g@addto@macro\th@definition{\thm@headpunct{}}
\g@addto@macro\th@plain{\thm@headpunct{}}
\g@addto@macro\th@remark{\thm@headpunct{}}
\makeatother
\usepackage{graphicx}
\usepackage{titlesec}
\titleformat{\section}{\normalsize\bfseries}{\hspace{0mm} \thesection.}{10 pt}{}
\usepackage{indentfirst}
\usepackage{setspace}
\usepackage{float}
\usepackage[round]{natbib}
\usepackage{url}
\usepackage{hyperref}
\hypersetup{pdfborder={0 0 0}}   

\providecommand{\Abstract}[1]{\noindent\textbf{Abstract} \\ #1}
\providecommand{\titlepaper}[1]{\begin{spacing}{1.15}\noindent\bf\Large#1\end{spacing}}
\providecommand{\keywords}[1]{\noindent\rule{0.3 \textwidth}{1 pt }\newline \textbf{Keywords:} #1}
\theoremstyle{definition}
\newtheorem{thm}{Theorem}

\numberwithin{equation}{section} 

\newenvironment{pf}
    {{\noindent  \textrm{\textbf{Proof:}}}~~}

\usepackage{mathrsfs} 	%%	Math. Character:	\mathscr{}
\usepackage{booktabs}

\DeclareMathOperator{\arcsec}{\arcsec\,}

\begin{document}
\thispagestyle{empty}
%% 	
%\logothailstat
%\vskip2cm
%%----------		title		----------%%
\titlepaper{%
A Compounded Probability Model for Decreasing Hazard and its  \\
Inferential Properties 
}	%%

%%----------		name of the authors		----------%%
\noindent
%%
%%%%%		for only one author 
%%%
%Name1 Surname\textsuperscript{*}
%%%
%\vskip 1pt \noindent%%
%%%
%Affiliation Dept/Program/Center, Institution Name, City, State, Country
%%%
%%%%%
%%
%%%%%		for two authors
%%
\fontsize{11pt}{11pt}\textbf{Brijesh P. Singh\textsuperscript{1}, 
Utpal Dhar Das\textsuperscript{1,*}} \textbf{and Sandeep Singh}\textsuperscript{1}
\vskip 1pt \noindent\textsuperscript{1}%%
Department of Statistics, Institute of Science, Banaras Hindu University, Varanasi-221105 INDIA
%%
%\vskip 1pt \noindent\textsuperscript{2}%%
%Affiliation Dept/Program/Center, Institution Name, City, State, Country
%%
%%%%%
%%
%%%%%		for three authors
%%%
%Name1 Surname,\textsuperscript{\hspace{-3 pt}1,*} 
%Name2 Surname\textsuperscript{2} and
%Name3 Surname\textsuperscript{3}
%%%
%\vskip 1pt \noindent\textsuperscript{1}%%
%Affiliation Dept/Program/Center, Institution Name, City, State, Country
%%%
%\vskip 1pt \noindent\textsuperscript{2}%%
%Affiliation Dept/Program/Center, Institution Name, City, State, Country
%%%
%\vskip 1pt \noindent\textsuperscript{3}%%
%Affiliation Dept/Program/Center, Institution Name, City, State, Country
%%%
%%%%%
\vskip 1.15pt 
\noindent
\textsuperscript{*}Corresponding author; Email: utpal.statmath@gmail.com

%%----------		abstract		----------%%
\Abstract{\hspace*{12 pt}	%%
There are some real life issues that are exists in nature which has early failure. This type of problems can be modelled either by a complex distribution having more than one parameter or by finite mixture of some distribution. In this article a single parameter continuous distribution is introduced to model such type of problems. The base line distribution is exponential and it is compounded by lindley distribution. Some important properties of the proposed distribution such as distribution function, survival function, hazard function and cumulative hazard function are derived. The maximum likelihood estimate of the parameter is obtained which is not in closed form, thus iteration procedure is used to obtain the estimate of parameter. The moments of the proposed distribution does not exist thus median and mode is obtained. The distribution is positively skewed and the hazard rate of this distribution is decreasing. Some real data sets are used to see the performance of proposed distribution with comparison of some other competent distributions of decreasing hazard using Likelihood, AIC, AICc, BIC and KS statistics.
}	%%

%%----------		keywords		----------%%
\vskip 5 pt

%%----------		contents		----------%%
%% 	section 1
\section{Introduction}
Lifetime distributions are used to describe statistically, length of the life of a system or a device. Normal, exponential, gamma and weibull distributions are the basic distributions that demonstrated in a number of theoretical results in the distributions theory. Particularly, exponential distribution is an invariable example for a number of theoretical concepts in reliability studies. It is characterized as constant hazard rate. In case of necessity for an increasing/decreasing failure rate model ordinarily the choice falls on weibull distribution. Lindley distribution is an increasing hazard rate distribution and has its own importance as a life testing distribution. The lindley distribution is one parameter distribution that is a mixture of exponential and gamma distributions and was proposed by \cite{lindley1958} in the context of Bayesian statistics. The lindley distribution is used to explain the lifetime phenomenon such as engineering, biology, medicine, ecology and finance. \cite{ghitany2011two} stated that it is mainly useful for modeling in mortality studies. Lindley distribution has generated little attention in excess of the exponential distribution because of its decreasing mean residual life function and increasing hazard rate however exponential distribution has constant mean residual life function and hazard rate.
A compounding of exponential and lindley distribution is considered in this study and named as compounded exponential-lindley (CEL) distribution. The distributional properties, estimation of parameters, Fisher informatio, entropies, stochastic ordering, quantile function, order statistics and simulation study for the proposed distribution have been discussed in detail. \\ 
	
\cite{adamidis1998lifetime} introduced a two-parameter lifetime distribution with decreasing hazard rate by compounding exponential and geometric distributions and \cite{barreto2011weibull} introduced a lifetime model with decreasing failure rate by compounding exponential and poisson-lindley distribution (EPL) with probability density function is given by
\begin{equation}
	f_{epl}(x;\beta,\theta)=\frac{\beta\theta^2(1+\theta)^2e^{-\beta x}}{(1+3\theta+\theta^2)}\frac{(3+\theta-e^{-\beta x})}{(1+\theta-e^{-\beta x})^3};x>0,\beta>0,\theta>0
\end{equation}
Another idea was proposed by \cite{kucs2007new} and \cite{tahmasbi2008two}. They introduced the exponential Poisson (EP) and exponential logarithmic (EL) distributions and the pdf is given by
\begin{align}
	f_{ep}(x;\beta,\lambda)=\frac{\lambda\beta}{1-e^{-\lambda}}e^{-\lambda-\beta x+\lambda e^{-\beta x}};x>0,\beta>0,\lambda>0 \\
	f_{el}(x;\beta,p)=\frac{1}{-\log{p}}\frac{\beta(1-p)e^{-\beta x}}{1-(1-p)e^{-\beta x}};x>0,\beta>0,p\in(0,1)
\end{align}
\cite{chahkandi2009some} introduced a class of distributions, which is exponential power series distributions (EPS), where compounding procedure follows the same way that was previously given by \cite{adamidis1998lifetime}.
\cite{weibull1951} a Swedish mathematician Waloddi Weibull describe the weibull distribution and the pdf is defined as
\begin{equation}
	f_{w}(x;\beta,\alpha)=\alpha\beta^{\alpha}x^{\alpha-1}e^{-\beta x};x>0,\beta>0,\alpha>0
\end{equation}
Natural mixing of exponential populations, giving rise to a decreasing hazard rate distribution, were first introduce by \cite{proschan1963theoretical}. Subsequently other distributions with decreasing hazard rates of practical interest were discussed by \cite{cozzolino1968probabilistic}.
The distributions with decreasing failure rate (DFR) are discussed in the works of \cite{lomax1954business}, \cite{barlow1963properties}, \cite{barlow1964bounds,barlow1965tables}, \cite{marshall1965maximumy},  \cite{dahiya1972goodness}, \cite{saunders1983maximum}, \cite{nassar1988two}, \cite{gleser1989gamma}, \cite{gurland1994shorter}. Keeping these ideas in view in this study, we introduce a new lifetime distribution by compounding exponential and lindley distribution.

\section{Proposed Distribution}
Let $X_1,X_2,...,X_n$  be a random sample from following exponential distribution with scale parameter $\lambda>0$ and the probability density function is in the form
\begin{align}
	f(x|\lambda)=\lambda e^{-\lambda x};\quad x>0,\lambda>0
\end{align}
The parameter $\lambda>0$ of the above distribution takes continuous value and hazard of the distribution is constant. Now we assume the parameter $\lambda$ is a random variable following lindley distribution. The probability density function of lindley distribution is given below
\begin{align}
	\phi(\lambda;\theta)=\frac{\theta^2}{(\theta+1)}(1+\lambda)e^{-\theta\lambda};\quad\theta>0,\lambda>0
\end{align}
Now the probability density function of the proposed distribution is given by 
\begin{align}\label{eq:pdf}
	g(x;\theta)&=\int\limits_{0}^{\infty}f(x|\lambda)\phi(\lambda;\theta)d\lambda=\frac{\theta^2}{(\theta+1)}\int\limits_{0}^{\infty}(\lambda+\lambda^2)e^{-\lambda(x+\theta)}d\lambda\nonumber\\
	&=\frac{\theta^2}{(\theta+1)}\frac{(x+\theta+2)}{(x+\theta)^3};\quad x>0,\theta>0
\end{align}
Since in the proposed distribution the parameter $\lambda$ of the exponential distribution follows lindley distribution, thus the new distribution is named as compounded exponentia lindley (CEL) distribution.\\

\noindent Now we denote the random variable $X$ follows a single parameter compounded exponentia lindley (CEL) distribution i.e
\begin{align}
	X\sim CEL(\theta) \qquad\implies g(x;\theta)=\frac{\theta^2}{(\theta+1)}\frac{(x+\theta+2)}{(x+\theta)^3};x>0,\theta>0
\end{align}
\begin{figure}[H]
	\centering
	\includegraphics[width=0.75\linewidth]{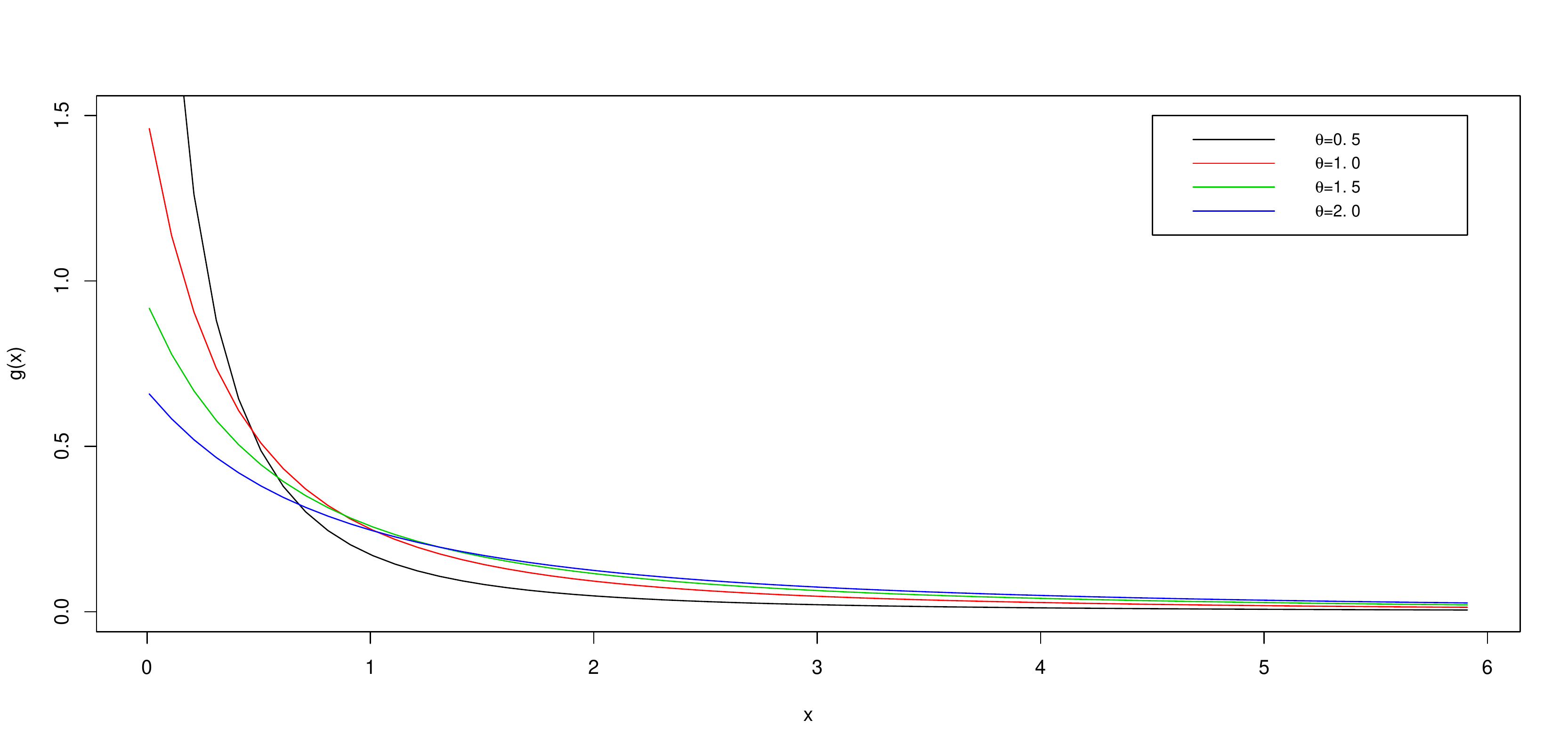}
	\caption{ Probability density function of CEL distribution}
	\label{fig:pdf}	
\end{figure}
and the cumulative distribution function (cdf) of CEL is obtained as
\begin{align}\label{eq:cdf}
	G(x;\theta)=\frac{x\left[x(\theta+1)+\theta(\theta+2)\right]}{(\theta+1)(x+\theta)^2};x>0,\theta>0
\end{align}
\begin{figure}[H]
	\centering
	\includegraphics[width=0.75\linewidth]{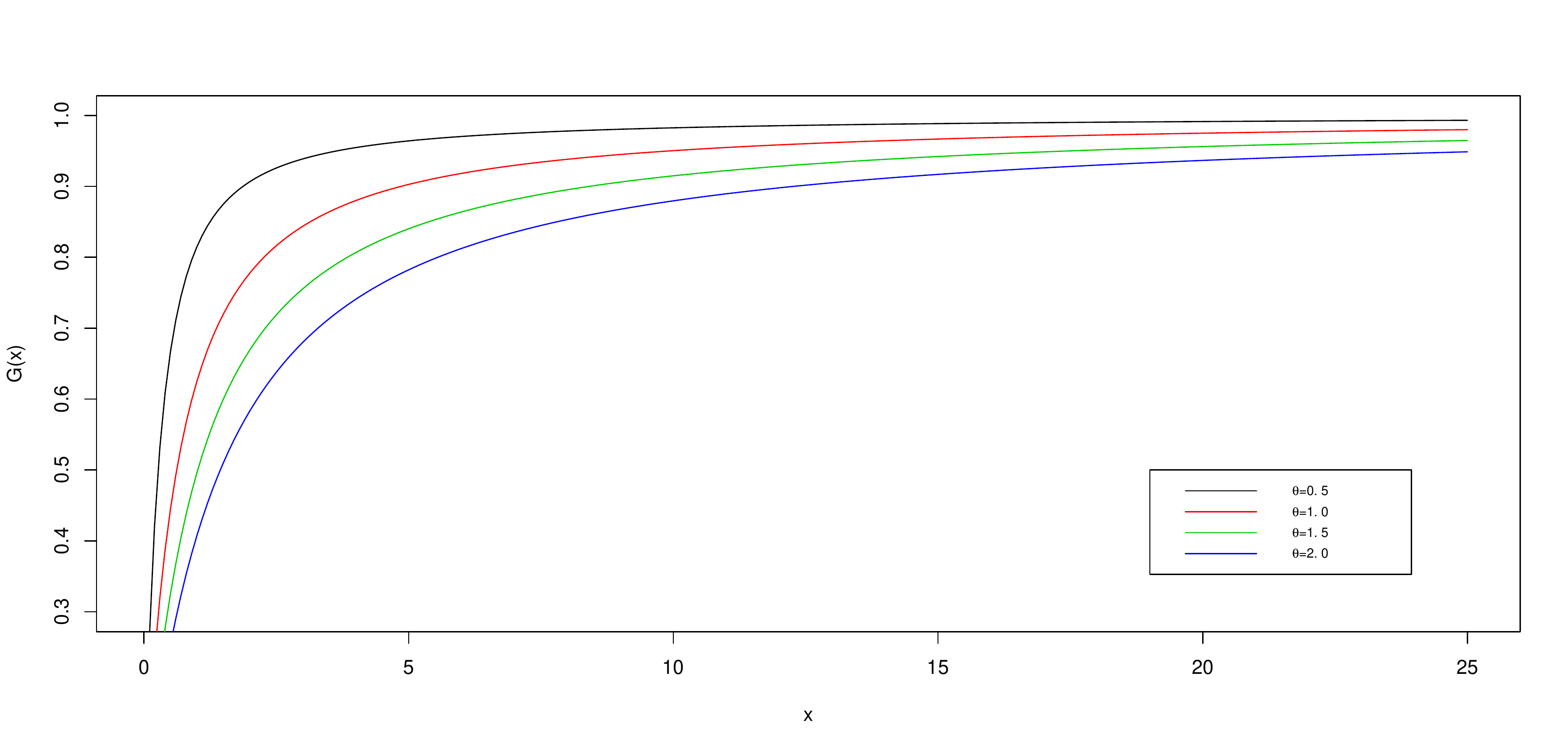}
	\caption{Cumulative distribution function of CEL distribution}
	\label{fig:cdf}	
\end{figure}
From the figure 1 and 2, it is clear that the distribution is early failure distribution for smaller value of $\theta$. The survival function $ S(x) $ of CEL having pdf \eqref{eq:pdf}, is given as
\begin{align}
	S(x)=\frac{\theta^2(x+\theta+1)}{(\theta+1)(x+\theta)^2}
\end{align}
\begin{figure}[H]
	\centering
	\includegraphics[width=0.75\linewidth]{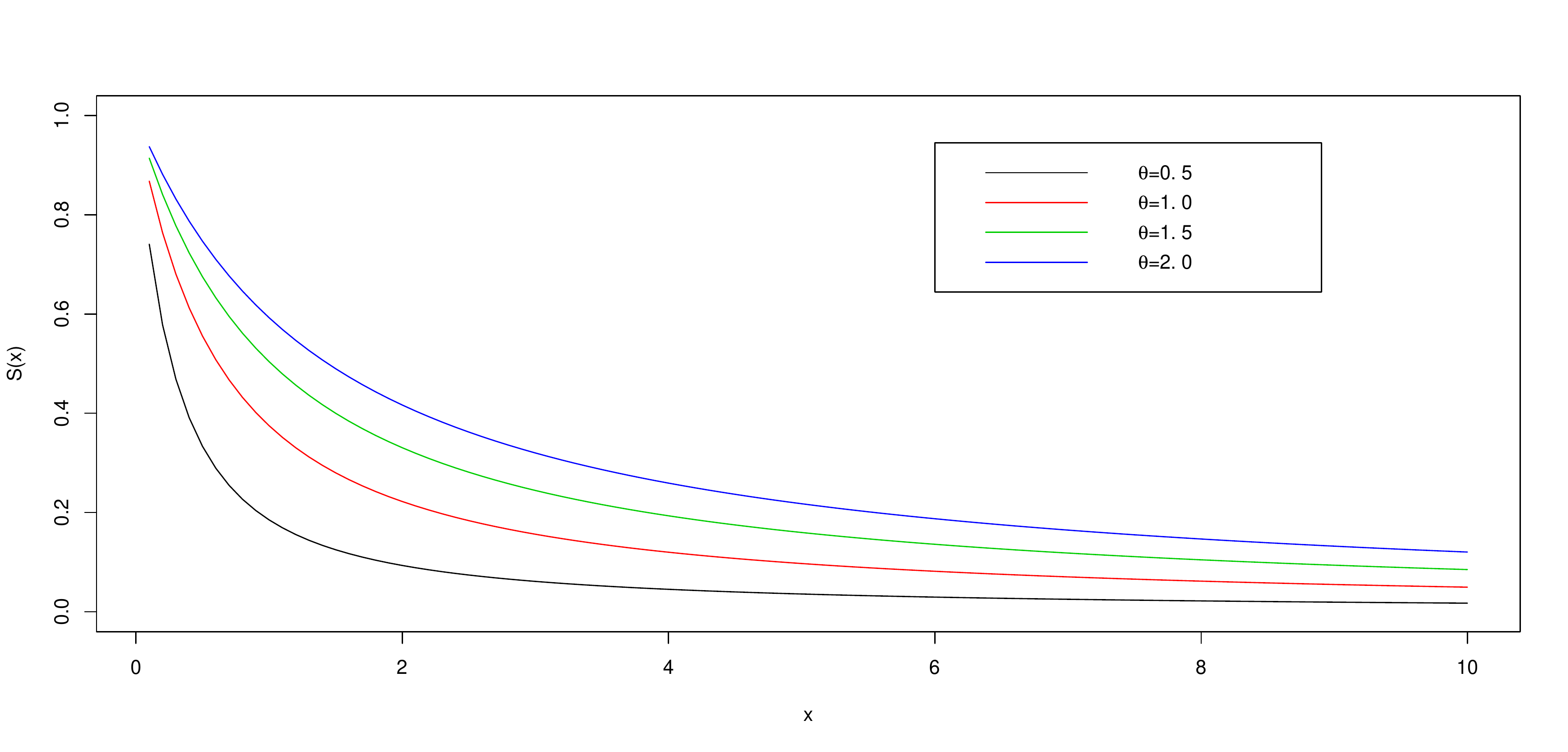}
	\caption{Survival function of CEL distribution}
	\label{fig:survival}	
\end{figure}
\noindent The hazard function is defined as \begin{align}\label{eq:haz}
	h(x)=\frac{g(x)}{S(x)}=\frac{g(x)}{1-G(x)}=\frac{(x+\theta+2)}{(x+\theta)(x+\theta+1)}
\end{align}
\begin{figure}[htbp]
	\centering
	\includegraphics[width=0.75\linewidth]{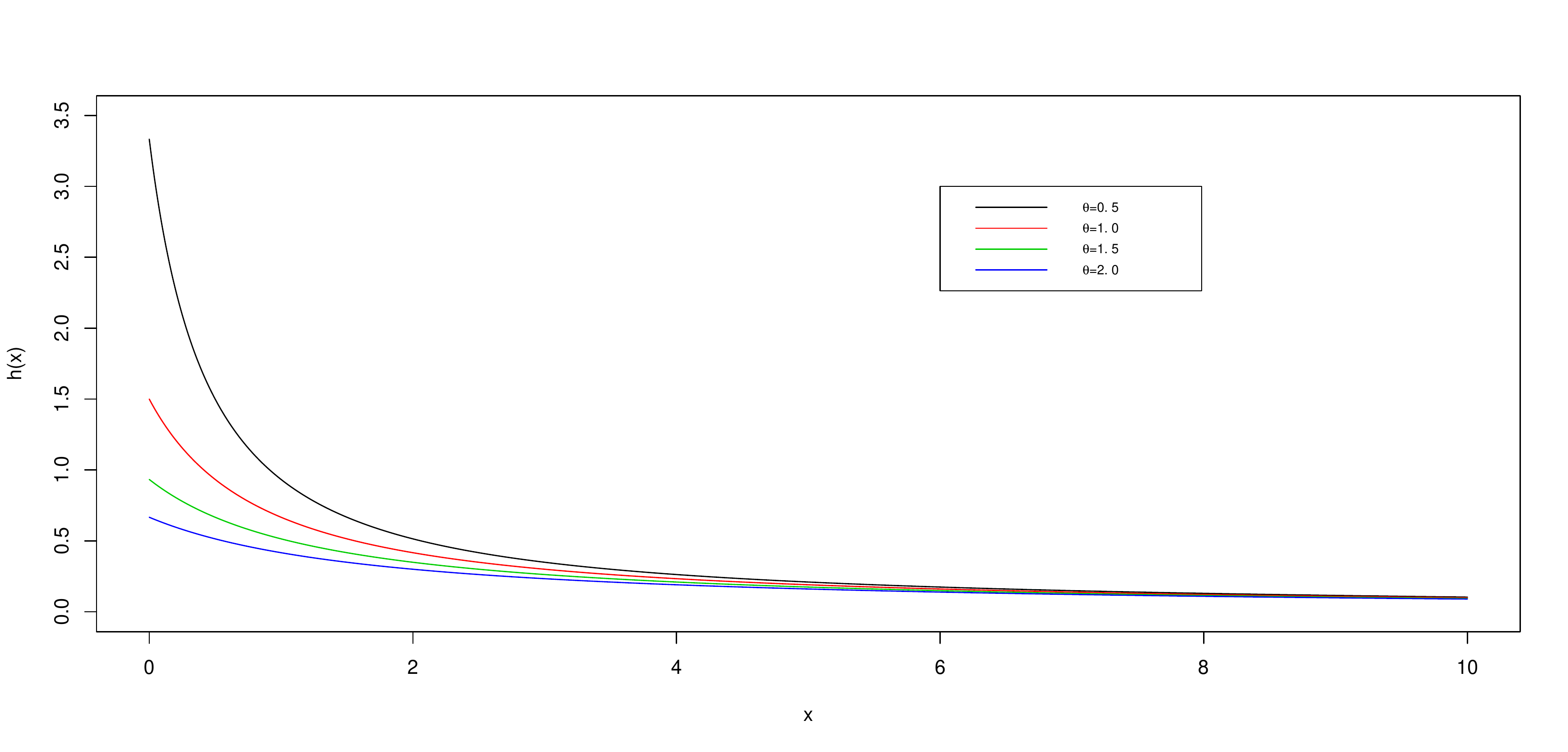}
	\caption{Hazard rate function of CEL distribution}
	\label{fig:hazard}	
\end{figure}
According to \cite{laser1980bathtub} we can find the shape of hazard rate, if $ g(t) $ is density function and $ g'(t) $ is the first order derivative of $ g(t) $ with respect to $ t $. He defined the term $ \eta(t)=-\frac{g'(t)}{g(t)} $, where $ g(t) $ is density function. Then he proved that if $ \eta'(t)>0 \quad\forall \quad t>0$, then the distribution has increasing failure rate (IFR) and if  $ \eta'(t)<0 \quad\forall \quad t>0$, then the distribution has decreasing failure rate (DFR).
\begin{align}
	\eta(t)=\frac{2t+3\theta+3}{(t+\theta)(t+\theta+2)}
\end{align}
Differentiating $ \eta(t) $ with respect to $ t $ we get
\begin{align}\label{eq:dfr}
	\eta'(t)=-\frac{2}{(t+\theta)(t+\theta+2)}-\frac{4}{(t+\theta)^2(t+\theta+2)}-\frac{2(\theta-1)}{(t+\theta)(t+\theta+2)^2}\nonumber\\-\frac{2(\theta-1)}{\left[(t+\theta)(t+\theta+2)\right]^2}
\end{align}
Now from the equation \eqref{eq:dfr} we have $ \eta'(t)<0 $ for all $ t>0 $, hence distribution has DFR.
\begin{align}
	h(x)=\frac{2}{(x+\theta)}-\frac{1}{(x+\theta+1)}\nonumber
\end{align}
After differentiating \eqref{eq:haz} with respect to $ x $ we get
\begin{align}
	h'(x)=-\frac{2}{(x+\theta)^2}+\frac{1}{(x+\theta+1)^{2}}\nonumber\\
	\lim_{x\to{0}}h'(x)=-\frac{2}{\theta^2}+\frac{1}{(\theta+1)^{2}}<0\quad \forall\quad \theta>0
\end{align}
Therefore $ h'(0)<0 \quad \forall\quad \theta>0$,
Hence CEL is a distribution of monotonic decreasing hazard with increasing time, however, the hazard function of exponential distribution is constant. CEL can be use for a phenomena having decreasing hazard such as infant mortality. The hazard function is shown in Figure~\ref{fig:hazard} indicates that, as the value of $\theta$  is increasing the hazard of the distribution becomes flatten.\newline
Hence, $ h(x) $  is a decreasing hazard rate function.

Now Cumulative hazard function $ H(t) $ is defined as
\begin{equation}
	H(t)=\int\limits_{0}^{t}h(x)dx=\log\left[\left(\frac{\theta+1}{t+\theta+1}\right)\left(\frac{t+\theta}{\theta}\right)^{2}\right]
\end{equation}
\begin{thm}
	The moments of the $CEL(\theta)$ distribution does not exists.
\end{thm}
\begin{pf}
	Suppose the random variable $ X $ comes from $CEL(\theta)$ then the $ r^{th} $ moment is given by
	\begin{align}
		E(X^r)=\int\limits_{0}^{\infty}x^rg(x)dx=\frac{\theta^2}{\theta+1}\int\limits_{0}^{\infty}x^r\frac{x+\theta+2}{(x+\theta)^3}dx\nonumber
	\end{align}
	Now
	\begin{align}
		\frac{1}{\theta+1}\int\limits_{0}^{\infty}\frac{x^r}{\left(1+\frac{x}{\theta}\right)^2}dx+\frac{2}{\theta(\theta+1)}\int\limits_{0}^{\infty}\frac{x^r}{\left(1+\frac{x}{\theta}\right)^3}dx\nonumber
	\end{align}
	Let $ \frac{x}{\theta}=z;\quad dx=\theta dz;\quad x\to0,z\to{0},\text{and}\quad x\to\infty,z\to{\infty} $ above integral become 
	\begin{align}
		\frac{\theta^{r+1}}{\theta+1}\int\limits_{0}^{\infty}\frac{z^r}{\left(1+z\right)^2}dz+\frac{2\theta^{r+1}}{\theta(\theta+1)}\int\limits_{0}^{\infty}\frac{z^r}{\left(1+z\right)^3}dz\nonumber
	\end{align}
	Now using Beta integral of second kind i.e $ \int\limits_{0}^{\infty}\frac{x^{m-1}}{\left(1+x\right)^{m+n}}dx=B(m,n)\quad ;m>0;n>0 $, we get 
	\begin{align}
		E(X^r)&=\frac{\theta^{r+1}}{\theta+1}B(r+1,1-r)+\frac{2\theta^{r+1}}{\theta(\theta+1)}B(r+1,2-r)\nonumber\\
		&=\frac{\theta^{r+1}}{\theta+1}\left[B(r+1,1-r)+\frac{2}{\theta}B(r+1,2-r)\right]
	\end{align}
	Here range is  $ -1<r<1  $. But range of $r$ should be $ r\geq1  $. Hence $ E(X^r) $ does not exists. Therefore mean, variance, SD as well as higher order moments does not edxists for $CEL(\theta)$.
\end{pf}

\begin{thm}
	The moment generating function of $CEL(\theta)$ does not exists.
\end{thm}
\begin{pf}
Let $ X $ be the random variable from $NWEL(\theta)$ distribution then the moment generating function (mgf) is given by
\begin{align}
	E(e^{tx})&=\int\limits_{0}^{\infty}e^{tx}g(x)dx=\frac{\theta^2}{\theta+1}\int\limits_{0}^{\infty}e^{tx}\frac{x+\theta+2}{(x+\theta)^3}dx\nonumber\\
	&=\frac{\theta^2}{\theta+1}\left[\int\limits_{0}^{\infty}\frac{e^{tx}}{(x+\theta)^2}dx+\int\limits_{0}^{\infty}\frac{2e^{tx}}{(x+\theta)^3}dx\right]
\end{align}
Now
\begin{align}
	\int\limits_{0}^{\infty}\frac{e^{tx}}{(x+\theta)^2}dx&=\left[\frac{e^{tx}}{-(x+\theta)}\right]^{\infty}_{0}+t\int\limits_{0}^{\infty}\frac{e^{tx}}{(x+\theta)}dx\nonumber\\
	&=\frac{1}{\theta}+\lim\limits_{\epsilon\to\infty}\left[t\int\limits_{0}^{\epsilon}\frac{e^{tx}}{(x+\theta)}dx\right]
\end{align}
Now applying L'Hospital rules we get 
\begin{align}
	\lim\limits_{x\to\infty}\frac{e^{tx}}{(x+\theta)}=\lim\limits_{x\to\infty}\frac{te^{tx}}{1}=\infty\nonumber
\end{align}
Hence integrand is divergent, as well as the function is not integrable over R we conclude thta $ E(e^{tx}) $ does not exists.\\
\end{pf}
The characteristic function of CEL distribution is defined as
\begin{align}\label{eq:chr}
	\Phi_x(t)=\int\limits_{0}^{\infty}e^{itx}g(x)dx=\frac{1}{\theta+1}\sum_{k=0}^{\infty}(-1)^{k}\frac{(k+1)!}{(it)^{k+1}}\left[1+\frac{2}{\theta}(k+2)\right]
\end{align}

\section{Entropies}
An entropy is a measure of randomness of any system. Entropy is an important property of probability distributions and it measures the uncertainty in a probability distribution. 
\subsection{Rényi Entropye}
An entropy is a measure of variation of the uncertainty, \cite{renyi1961} gave an expression of the Entropy function defined by 
\begin{align}
	e(\eta)&=\frac{1}{1-\eta}\log\left[\int\limits_0^\infty g^\eta (x)dx\right]\nonumber
\end{align}
where $ 0<\eta<1 $, Substituting the value of g(x) from \eqref{eq:pdf}
\begin{align}
	e(\eta)&=\frac{1}{1-\eta}\log\left[\int\limits_0^\infty \left(\frac{\theta^2}{(\theta+1)}\frac{(x+\theta+2)}{(x+\theta)^3}\right)^\eta dx\right]\nonumber\\
	&=\frac{1}{1-\eta}\log\left[\left(\frac{\theta^2}{\theta+1}\right)^\eta\int\limits_0^\infty\left\{\frac{1}{(x+\theta)^2}+\frac{2}{(x+\theta)^3}\right\}dx\right]\nonumber
\end{align}
Now applying Binomial expansion $ (a+b)^n=\sum\limits_{k=0}^{n}{{n}\choose{k}}a^kb^{n-k} $ we get
\begin{align}
	\frac{1}{1-\eta}\log\left[\left(\frac{\theta^2}{\theta+1}\right)^\eta\int\limits_0^\infty\sum\limits_{k=0}^{\eta}{{\eta}\choose{k}}\left(\frac{1}{x+\theta}\right)^{2k}\left(\frac{2}{(x+\theta)^3}\right)^{\eta-k}dx\right]\nonumber
\end{align}
after simlification we get the Renyi entropy as
\begin{align}
	e(\eta)=\frac{\eta}{1-\eta}\log\left(\frac{\theta^2}{\theta+1}\right)+\frac{1}{1-\eta}\log\left[\sum\limits_{k=0}^{\eta}{{\eta}\choose{k}}\frac{2^{\eta-k}}{(3\eta-k-1)\theta^{(3\eta-k-1)}}\right]
\end{align}
where $ 0<\eta<1, \quad \theta>0,\quad x>0 $
\subsection{Tsallis Entropy}
This is introduced by \cite{tsallis1988possible} as a basis for generalizing the standard statistical mechanics
\begin{align}
	S_{\lambda}&=\frac{1}{1-\lambda}\left[1-\int\limits_0^\infty g^\lambda (x)dx\right]\nonumber\\
	&=\frac{1}{1-\lambda}\left[1-\left(\frac{\theta^2}{(\theta+1)}\right)^\lambda\int\limits_0^\infty \left(\frac{(x+\theta+2)}{(x+\theta)^3}\right)^\lambda dx\right]\nonumber
\end{align}
Now applying Binomial expansion $ (a+b)^n=\sum\limits_{k=0}^{n}{{n}\choose{k}}a^kb^{n-k} $ and simplifying we get Tsallis Entropy as in \eqref{eq:tsal}.
\begin{align}\label{eq:tsal}
	e(\eta)=\frac{1}{1-\lambda}\left[1-\left(\frac{\theta^2}{\theta+1}\right)^\lambda\sum\limits_{k=0}^{\lambda}{{\lambda}\choose{k}}\frac{2^{\lambda-k}}{(3\lambda-k-1)\theta^{(3\lambda-k-1)}}\right]
\end{align}
\section{Quantile Function}
The quantile function for CEL distribution is defined in the form $ x_q=Q(u)=G^{-1}(u) $ where $ Q(u) $ is the quantile function of $ G(x) $ in the range $ 0<u<1 $.\\
Taking $ G(x) $ is the cdf of CEL and inverting it as above will give us the quantile function as follows
\begin{align}\label{eq:quan}
	G(x)=\frac{x\left[x(\theta+1)+\theta(\theta+2)\right]}{(\theta+1)(x+\theta)^2}=u
\end{align}
Simplifying equation \eqref{eq:quan} above gives the following:
\begin{equation}
	\left(\frac{x}{x+\theta}\right)^2+\frac{x\theta(\theta+2)}{(x+\theta)^2}=u\nonumber
\end{equation}
Now let $ \frac{x}{x+\theta}=z $ we get from above 
\begin{align}
	z^2+\left(\frac{\theta+2}{\theta+1}\right)z(1-z)=u\nonumber\\
	z^2-z(\theta+2)+u(\theta+1)=0
\end{align}
This is a quadratic equation and after solving we get the solution for $ x $ as
\begin{align}\label{eq:quanf}
	z=\frac{x}{x+\theta}=\frac{(\theta+2)\pm\sqrt{(\theta+2)^2-4u(\theta+1)}}{2}\nonumber\\
	Q(u)=\theta\left[\frac{2}{-\theta\pm\sqrt{(\theta+2)^2-4u(\theta+1)}}-1\right]
\end{align}
where $ u $ is a uniform variate on the unit interval (0,1).\\
The median of $ X $ from the CEL is simply obtained by setting $ u=0.5 $ and this substitution of $ u = 0.5 $ in the above equation \eqref{eq:quanf} gives.
\begin{align}
	Median=\theta\left[\frac{2}{-\theta+\sqrt{(\theta+1)^2+1}}-1\right]
\end{align}
Bowley’s measure of skewness based on quartiles is defined as:
\begin{align}
	SK=\frac{Q(\frac{3}{4})-2Q(\frac{1}{2})+Q(\frac{1}{4})}{Q(\frac{3}{4})-Q(\frac{1}{4})}
\end{align}
and \cite{moors1988quantile} presented the Moors’ kurtosis based on octiles by
\begin{align}
	KT=\frac{Q(\frac{7}{8})-Q(\frac{5}{8})-Q(\frac{3}{8})+Q(\frac{1}{8})}{Q(\frac{6}{8})-Q(\frac{1}{8})}
\end{align}
where Q(.) is calculated by using the quantile function from equation \eqref{eq:quanf}.
\section{Stochastic Orderings}
Stochastic ordering of a continuous random variable is an important tool to judging their comparative behaviour. A random variable X is said to be smaller than a random variable Y.\\
(i) Stochastic order $X\leq_{st}Y$ if $F_{X}{(x)}\geq F_{Y}{(x)}$ for all x.\\
(ii) Hazard rate order $X\leq_{hr}Y$ if $h_{X}{(x)}\geq h_{Y}{(x)}$ for all x.\\
(iii) Mean residual life order $X\leq_{mrl}Y$ if $m_{X}{(x)}\geq m_{Y}{(x)}$ for all x.\\
(iv) Likelihood ratio order $X\leq_{lr}Y$ if $\frac{f_{X}{(x)}}{f_{Y}{(x)}}$ decreases in x.\\
The following results by \cite{shakedjg} are well known for introducing stochastic ordering of distributions
\begin{align}
	X\leq_{lr}Y\implies X\leq_{hr}Y\implies X\leq_{mrl}Y\nonumber\\
	i.e\qquad X\leq_{st}Y\nonumber
\end{align}
with the help of following theorem we claim that CEL distribution is ordered with respect to strongest likelihood ratio ordering

\begin{thm}
	Let $X\sim CEL(\theta_1)$ distribution and $Y\sim CEL(\theta_2)$ distribution. If $\theta_1>\theta_2 $ then $X\leq_{lr}Y$ and therefore $X\leq_{hr}Y$, $X\leq_{mrl}Y$ and $X\leq_{st}Y$.
\end{thm}
\begin{pf}
	We have 
	\begin{align}
		\frac{f_{X}(x)}{f_{Y}(x)}=\frac{\theta^2_1(\theta_2+1)}{\theta^2_2(\theta_1+1)}\left(\frac{x+\theta_1+2}{x+\theta_2+2}\right)\left(\frac{x+\theta_2}{x+\theta_1}\right)^{3};\qquad x>0\nonumber
	\end{align}
	Now taking log both side we get
	\begin{align}
		\log\left[\frac{f_{X}(x)}{f_{Y}(x)}\right]=\log\left[\frac{\theta^2_1(\theta_2+1)}{\theta^2_2(\theta_1+1)}\right]+\log\left(\frac{x+\theta_1+2}{x+\theta_2+2}\right)+3\log\left(\frac{x+\theta_2}{x+\theta_1}\right)\nonumber
	\end{align}
	By differentiating both side we get
	\begin{align}
		\frac{d}{dx}\log\frac{f_{X}(x)}{f_{Y}(x)}=\frac{\theta_2-\theta_1}{(2+\theta_1+x)(2+\theta_2+x)}+\frac{3(\theta_2-\theta_1)}{(x+\theta_1)(x+\theta_2)}\nonumber
	\end{align}
	Thus for $\theta_1>\theta_2 ,\frac{d}{dx}\log\frac{f_{X}(x)}{f_{Y}(x)}<0. $This means that $X\leq_{lr}Y$ and hence $X\leq_{hr}Y$, $X\leq_{mrl}Y$ and $X\leq_{st}Y$.
\end{pf}
\section{Distribution of order statistics}
Let $X_{1},X_{2},...,X_{m}$ be a random sample of size $m$ from CEL and let
$ X_{1;m}\leq X_{2;m}\leq...\leq X_{m;m} $ represent the corresponding order statistics. The pdf of $ X_{m;m} $ i.e $ r^{th} $ order statistics is given by
\begin{align}\label{eq:or1}
	g_{(r:m)}(x)=\frac{m!}{(r-1)!(m-r)!}G^{r-1}(x)\left[1-G(x)\right]^{m-r}g(x)\nonumber\\
	=Z\sum_{l=0}^{m-r}{{m-r}\choose_{l}}(-1)^lG^{r+l-1}(x)g(x)
\end{align}
where $ Z=\frac{m!}{(r-1)!(m-r)!} $ and $ g(x) $ and $ G(x) $ are pdf and cdf of CEL defined in \eqref{eq:pdf} and \eqref{eq:cdf} respectively.\\
Substituting for $ G(x) $  and $ g(x) $ in \eqref{eq:or1} and applying the general binomial expansion, we have
\begin{align}
	g_{(r:m)}(x)&=Z\sum_{l=0}^{m-r}{{m-r}\choose_{l}}(-1)^l\left[\frac{x\left[x(\theta+1)+\theta(\theta+2)\right]}{(\theta+1)(x+\theta)^2}\right]^{r+l-1}\frac{\theta^2}{(\theta+1)}\frac{(x+\theta+2)}{(x+\theta)^3}\nonumber\\
	&=Z\sum_{l=0}^{m-r}\sum_{k=0}^{(r+l-1)}{{m-r}\choose_{l}}{{r+l-1}\choose_{k}}C_{l;k}\frac{x^{2r+2l-k-2}(x+\theta+2)}{(x+\theta)^{2r+2l+1}}
\end{align}
where $ C_{l;k}=(-1)^{l}\left(\frac{\theta^2}{(\theta+1)}\right)^{k+1}\left(\frac{\theta+2}{\theta}\right)^k $.\\
Hence, the pdf of the minimum order statistic $ X_{(1)} $ and maximum order statistic $ X_{(n)} $ of the CEL are respectively given by,
respectively given by
\begin{align}
	g_{(1:m)}(x)=Z\sum_{l=0}^{m-1}\sum_{k=0}^{l}{{m-1}\choose_{l}}{{l}\choose_{k}}C_{l;k}\frac{x^{2l-k}(x+\theta+2)}{(x+\theta)^{2l+3}}
\end{align}
\section{Estimation of the Parameter of CEL}
For the estimation of the parameter we have used the most popular estimation procedure i.e. Maximum Likelihood Estimation procedure.\\ Suppose  $ X = (X_1,X_2,X_3,...,X_n) $ be an independently and identically distributed (iid) random variables of size $ n $ with pdf \eqref{eq:pdf} from CEL($\theta$). Then, the likelihood function based on observed sample $ X = (x_1,x_2,x_3,...,x_n) $ is defined as 
\begin{align}\label{eq:mle}
	L(\theta;x)=\left(\frac{\theta^2}{\theta+1}\right)^n\displaystyle\prod_{i=0}^{n}\frac{x_i+\theta+2}{(x_i+\theta)^3}
\end{align}
The log-likelihood function corresponding to \eqref{eq:mle} is given by
\begin{align}\label{eq:mle1}
	\log{L}=2n\log\theta-n\log(\theta+1)+\sum\limits_{i=0}^{n}\left\{\log(x_i+\theta+2)-3\log(x_i+\theta)\right\}
\end{align}
Now after differentiating partially \eqref{eq:mle1} with respect to $ \theta $ and equating zero, we get 
\begin{align}
	\frac{\partial(\log{L})}{\partial\theta}=0\nonumber
\end{align}
Hence, the log-likelihood equation for estimating $\theta$ is given by \eqref{eq:mle2},
\begin{align}\label{eq:mle2}
	\frac{2n}{\theta}-\frac{n}{(\theta+1)}+\sum\limits_{i=0}^{n}\left\{\frac{1}{(x_i+\theta+2)}-\frac{3}{(x_i+\theta)}\right\}=0
\end{align}
Above equation is not solvable analytically for $\theta$. Thus numerical iteration technique is used to get its numerical solution.\\
Fisher Information matrix can be estimated by
\begin{align}
	I(\hat{\theta})&=\left[\frac{-\partial^2}{\partial\theta^2}\log{L}\right]_{\theta=\hat{\theta}}\nonumber\\
	\frac{\partial^2}{\partial\theta^2}\log{L}&=-\frac{2n}{\theta^2}+\frac{n}{(\theta+1)^2}+\sum\limits_{i=0}^{n}\left\{\frac{3}{(x_i+\theta)^2}-\frac{1}{(x_i+\theta+2)^2}\right\}
\end{align}
For large samples, we can obtain the confidence intervals based on Fisher information matrix $I^{-1}(\hat{\theta})$ which provides the estimated asymptotic
variance for the parameter $\theta$. Thus, a two-sided $100(1-\alpha)\%$  confidence interval of $\theta$ and it is defined as $\hat{\theta}\pm Z_{\alpha}/2\sqrt{var{\hat{\theta}}}$. Where $Z_{\alpha}/2$ denotes the upper $\alpha$-th percentile of the standard normal distribution.
\section{Simulation study}
In this section we evaluate the performance of the MLEs of the model parameter for the CEL distribution. We generate random variables from CEL($\theta$) and
then obtain m.l.e. of the parameter $\theta$, Now for $\theta$=2 we generate the sample size 20, 30, 50, 90, 150, 200. The program is replicated N= 2,500 times to get the maximum likelihood estimate of $\theta$. The simulation results are reported in Table~(\ref{table:tab1}).
\begin{table}[H]
	\centering
	\caption{Simulation results for $\theta$=2}\label{tb:data1}
	\tabcolsep=20pt
	\begin{tabular}{ccccc}
		\hline
		n & Bias & MSE & Var. & Est.\\\hline
		20 & 0.11554 & 0.53363 & 0.50319 & 2.11554\\
		30 & 0.07354 & 0.29175 & 0.30196 & 2.07354\\
		50 & 0.04340 & 0.16161 & 0.16965 & 2.04340\\
		90 & 0.01612 & 0.08701 & 0.08982 & 2.01612\\
		150& 0.01415 & 0.05453 & 0.05321 & 2.01415\\
		200& 0.00896 & 0.03761 & 0.03949 & 2.00896
		\\\hline
		\label{table:tab1}
	\end{tabular}
\end{table}
 We clearly observe from the Table~(\ref{table:tab1}) that the values of the bias and the mean square error (MSE) of the parameter estimates decreases as the sample size n increases. the MSE decreases as sample size n increases, it proves the consistency of the estimator.
\section{Goodness of fit}
The application of goodness of fit of CEL distribution has been discussed with two real data sets. First data set presents the results of a life-test experiment in which specimens of a type of electrical insulating fluid were subject to a constant voltage stress(34 KV/minutes), this data set is reported by \cite{nelson1982weibull} and other data is represents 30 failure times of the air conditioning system of an airplane has been reported in a paper by \cite{linhart1986} and has also been analyzed by \cite{barreto2013new} and so on.
For comparing the suitability of the model, we have considered following criterion’s; namely AIC (Akaike Information Criterion), BIC (Bayesian information criterion), AICc (Corrected Akaike information criterion)  and KS statistics with associated $p$-value of the fitted distributions are presented in Table~(\ref{table:tab2}) and Table~(\ref{table:tab3}).The AIC, BIC, AICc and KS Statistics are computed using the following formulae
\begin{align}
	AIC&=-2loglik+2k,\qquad BIC=-2loglik+k\log n\nonumber\\
	AICc&=AIC+\frac{2k^2+2k}{n-k-1},\qquad  D=\sup\limits_{x}|F_{n}{(x)}-F_{0}{(x)}|\nonumber
\end{align}
where $k$= the number of parameters, $n$= the sample size, and the $F_{n}{(x)}$=empirical distribution function and $ F_{0}{(x)} $   is the theoretical cumulativedistribution function.

\begin{table}[H]
	\centering
	\caption{MLE's, - 2ln L, AIC, KS and $p$-values of the fitted distributions for the 1st dataset.}
	\tabcolsep=7.5pt
	\begin{tabular}{cccccccc}
		\hline
		Distribution & Estimate & -2LL & AIC & BIC & AICc & KS & $p$-value
		\\\hline
		CEL($\theta$)  & 7.0385  & 137.98 & 139.98 & 140.92 & 140.21 & 0.1131 &  0.9458\\
		EPL($\beta,\theta$)  & (0.0334, 0.5521) & 136.18  & 140.18 & 142.06 & 140.93 & 0.1500 &    0.7312\\
		EL($\beta,p$)  & (0.0393, 0.0982)  & 135.98 & 139.98 & 141.87 & 140.73 & 0.1382 &   0.8137\\
		EP($\beta,\lambda$)  & (0.0409, 2.2112) & 136.89  & 140.89 & 142.78 & 141.64 & 0.1611 &   0.6497\\
		Weibull($\beta,\theta$)  & (0.0818, 0.7708) & 136.77  & 140.77 & 142.66 & 141.52 & 0.1613 &   0.6482\\
		Gamma($\beta,\theta$)  & (0.0480, 0.6897) & 137.23 & 141.23  & 143.12 & 141.98 & 0.1846 & 0.4802\\
		\hline
		\label{table:tab2}
	\end{tabular}
\end{table}
\begin{figure}[H]
	\centering
	\includegraphics[width=0.85\linewidth]{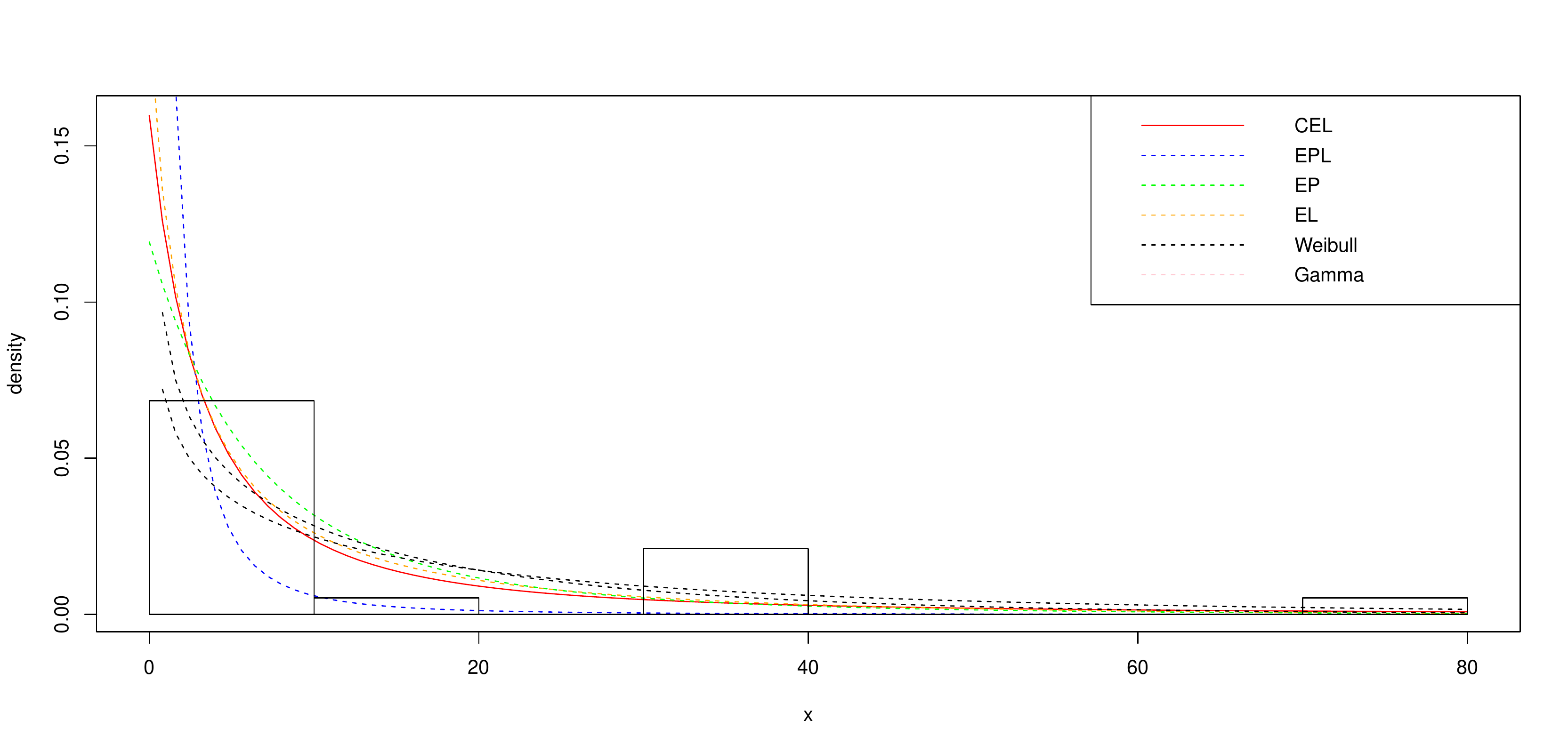}
	\caption{Fitted pdfs of 1st data set}
	\label{fig:1stp}	
\end{figure}
\begin{figure}[H]
	\centering
	\includegraphics[width=0.80\linewidth]{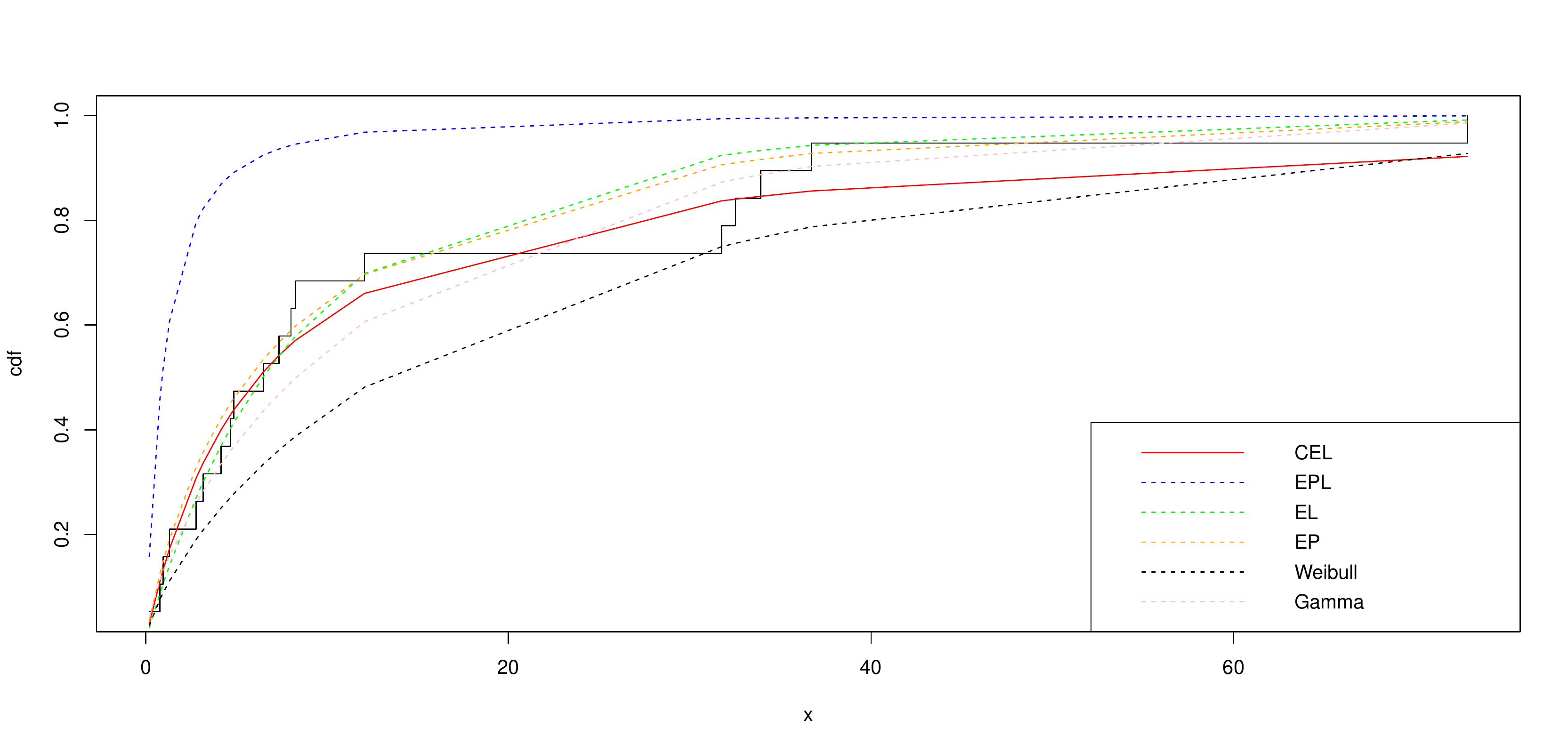}
	\caption{Fitted cdfs and ecdf of 1st data set}
	\label{fig:1stc}	
\end{figure}
\begin{figure}[H]
	\centering
	\includegraphics[width=0.95\linewidth]{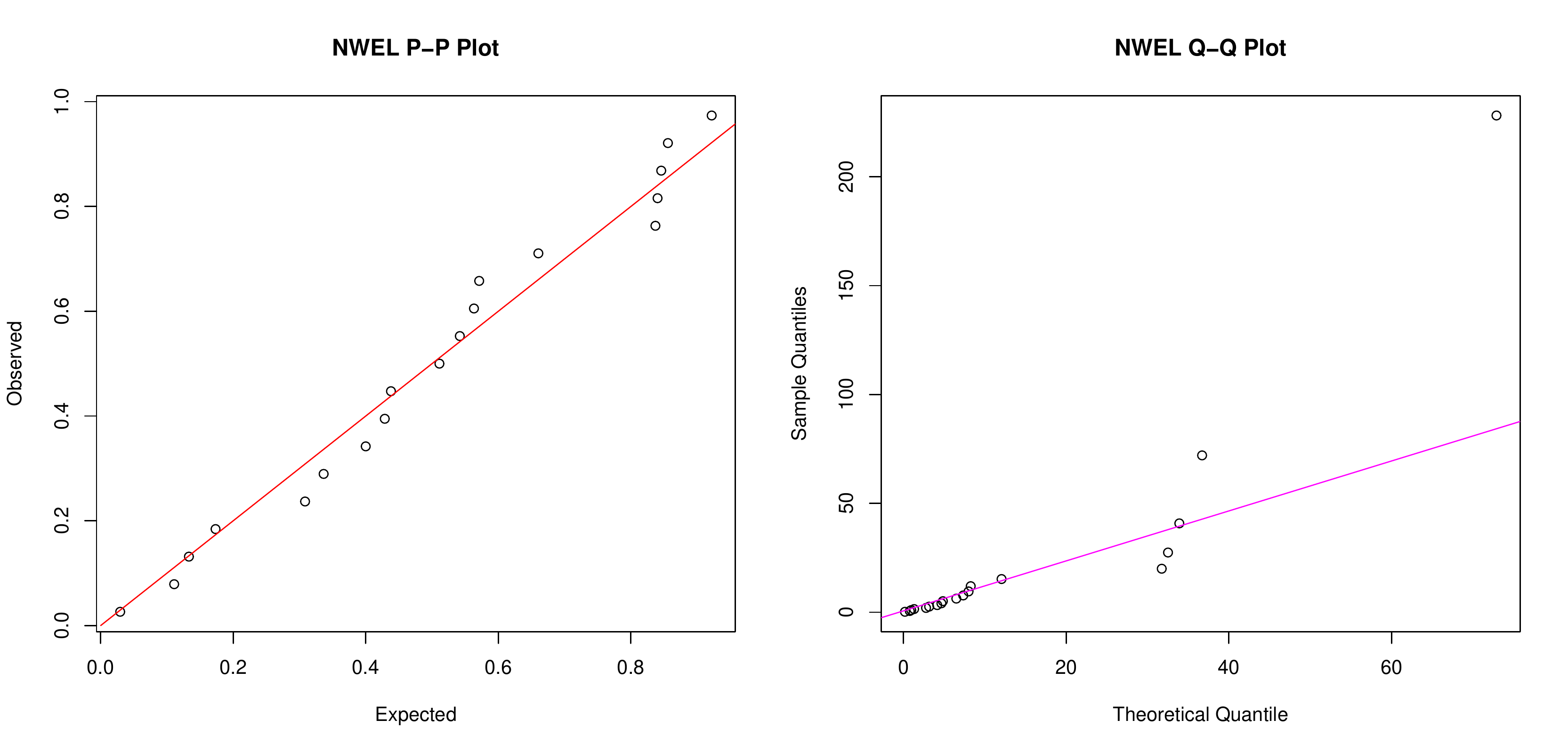}
	\caption{p-p and q-q plot for the 1st data set.}
	\label{fig:pp1}	
\end{figure}
Here we notice that all of the considered models fit the data at 5\% level of significance. the proposed distribution has minimum K-S and higher associated p-value among all the fitted models. We may say that proposed CEL distribution is the most acceptable model for the present data set according to KS statistics. For better visualization of the fitted models the estimated pdfs, cdfs, pp-qq plots are shown in Figure~\ref{fig:1stp}, Figure~\ref{fig:1stc}, Figure~\ref{fig:pp1} for the first data set. 
\begin{table}[H]
	\centering
	\caption{MLE's, - 2ln L, AIC, KS and $p$-values of the fitted distributions for the 2nd dataset.}\label{tb:data3}
	\tabcolsep=7.5pt
	\begin{tabular}{cccccccc}
		\hline 
		Distribution & Estimate & -2LL & AIC & BIC & AICc & KS & $p$-value
		\\\hline
		CEL($\theta$)  & 30.267  & 307.17 & 309.17 & 310.57 & 309.31 & 0.1061 &  0.8695\\
		EPL($\beta,\theta$)  & (0.0101, 0.9193) & 302.87  & 306.87 & 309.68 & 307.32 & 0.1282 &   0.7076\\
		EL($\beta,p$)  & (0.0111, 0.1932)  & 302.83 & 306.83 & 309.63 & 307.28 & 0.1291 &   0.6986\\
		EP($\beta,\theta$)  & (0.0105, 1.8243) & 303.22  & 307.22 & 310.02 & 307.66 & 0.1468 &   0.5375\\
		Weibull($\beta,\theta$)  & (0.0183, 0.8536) & 307.87  & 310.68 & 308.32 & 303.87 & 0.1534 &   0.4806\\
		Gamma($\beta,\theta$)  & (0.0136, 0.8119) & 304.33 & 308.33  & 311.13 & 308.78 & 0.1694 & 0.3556\\
		\hline
		\label{table:tab3}
	\end{tabular}
\end{table}
\begin{figure}[H]
	\centering
	\includegraphics[width=0.85\linewidth]{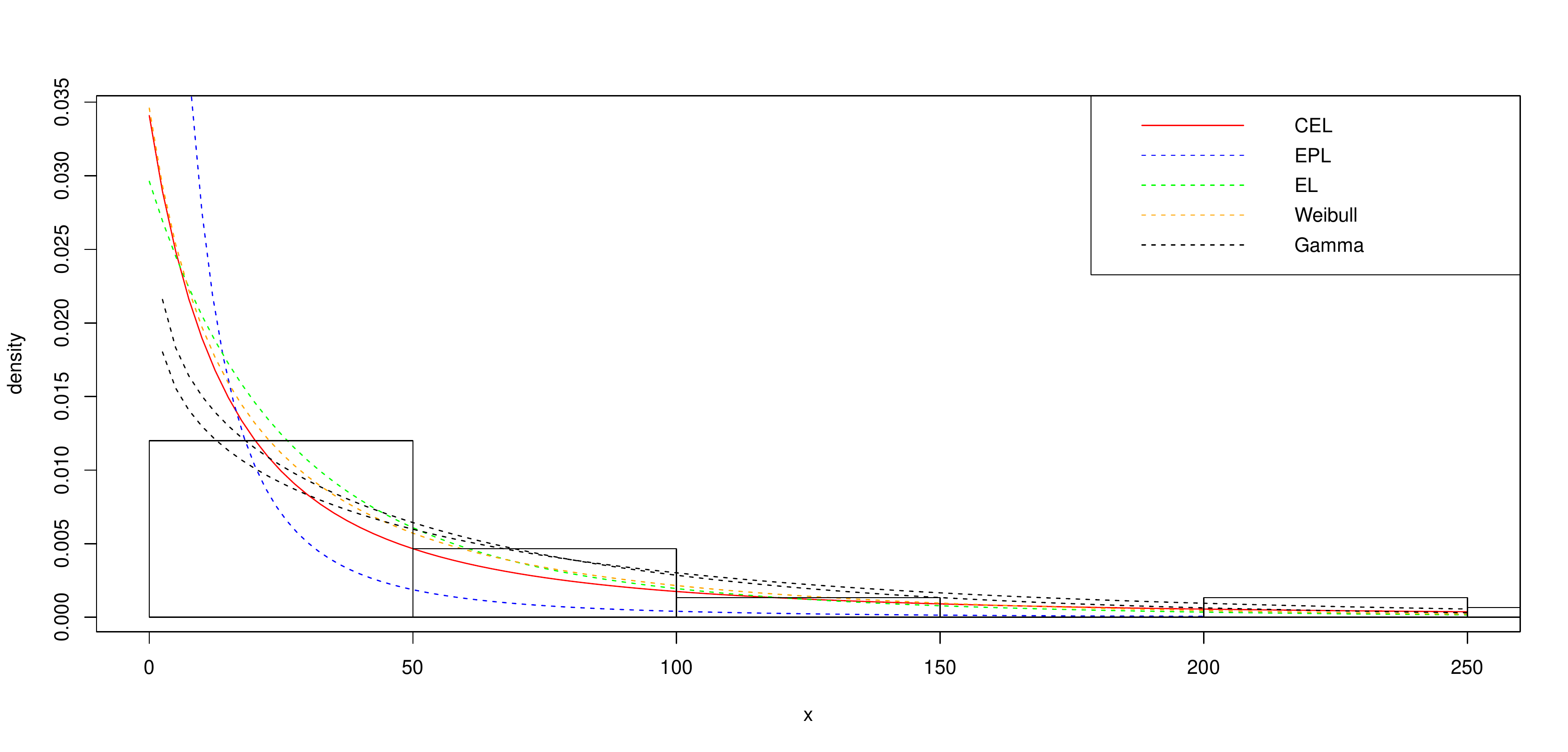}
	\caption{Fitted pdfs of 2nd data set}
	\label{fig:2ndp}	
\end{figure}
\begin{figure}[H]
	\centering
	\includegraphics[width=0.85\linewidth]{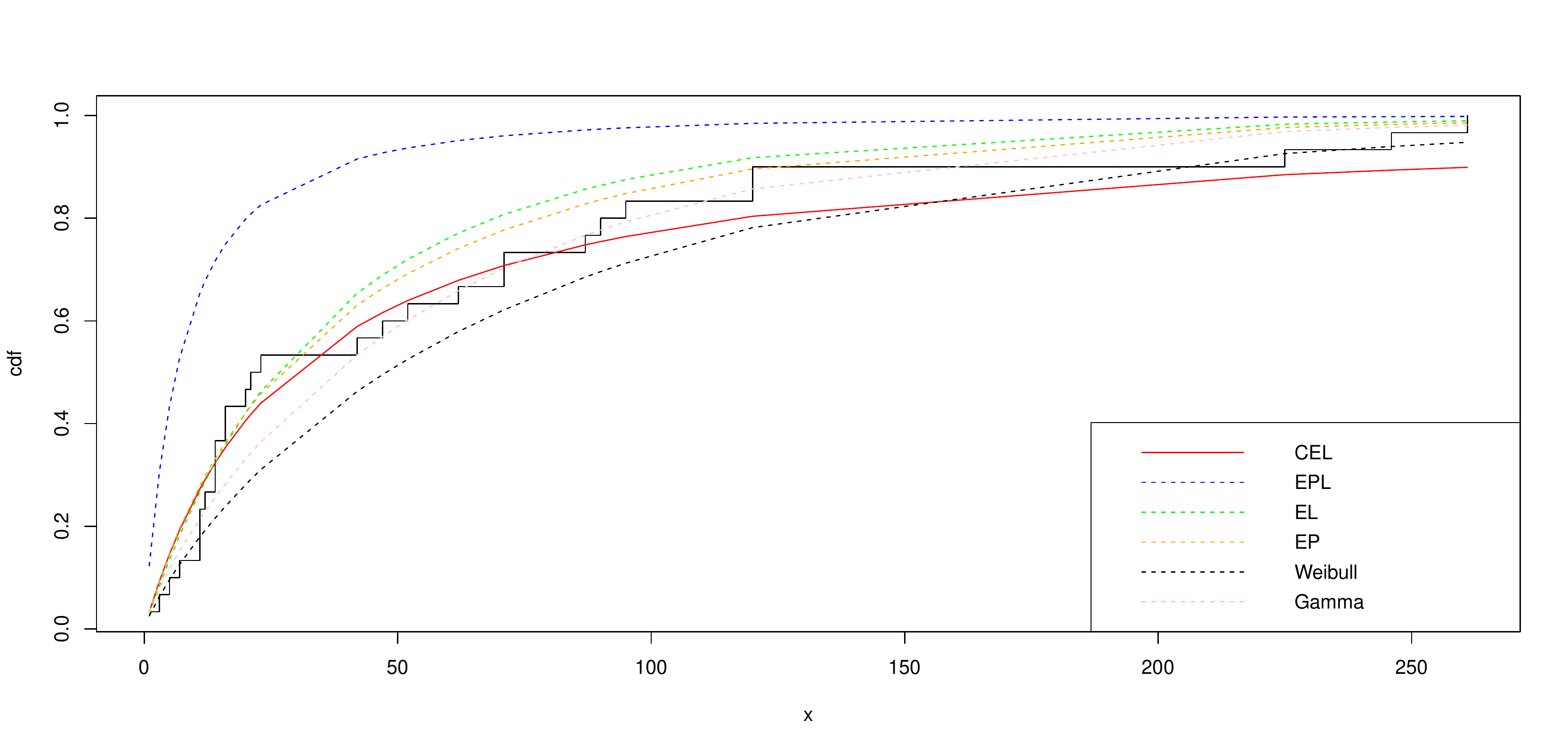}
	\caption{Fitted cdfs of 2nd data set}
	\label{fig:2ndc}	
\end{figure}
\begin{figure}[H]
	\centering
	\includegraphics[width=0.85\linewidth]{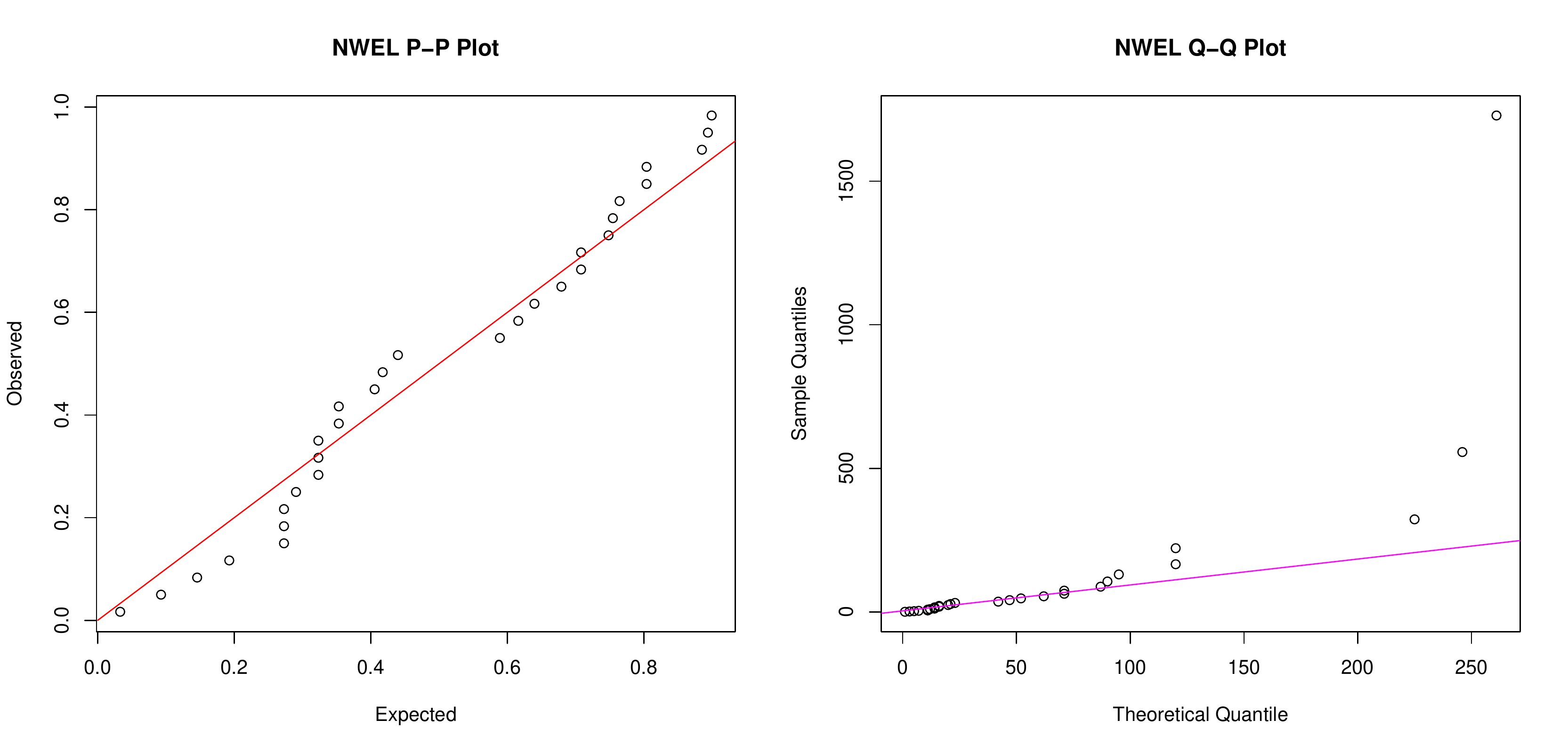}
	\caption{p-p and q-q plot for the 2nd data set.}
	\label{fig:pp2}	
\end{figure}
For the second data set also all the considered models fit well. Here also, the value of KS statistics is the minimum for CEL distribution with the maximum associated $p$-value. From the above discussion on two real data sets we see that all the considered six decreasing failure models fit to the two datasets considered here. The proposed distribution has minimum KS statistics with associated higher p-value. However for the model selection criterion  i.e. AIC, BIC, AICc for both data sets the Exponential Logarithmic (EL) distribution is performed well. Moreover, all the considered two-parameter model whereas Proposed model has a single parameter distribution. For better visualization of the fitted models the estimated pdfs, cdfs, pp-qq plots are shown in Figure~\ref{fig:2ndp}, Figure~\ref{fig:2ndc}, Figure~\ref{fig:pp2} for the second data set.These figures also indicate the proposed CEL distribution is the better fit to the considered real data in comparison to other competitive decreasing failure rate models.

\section{Conclusions}
A single parameter lifetime distribution CEL($\theta$) distribution has been introduced by weighted Lindley distribution. The CEL($\theta$) distribution is mean free distribution and single parameter decreasing hazard. The moment generating function, $r^{th}$ oreder moments does not exists thus mean,  variance, cumulant generating function, mean deviation about mean and median, Bonferroni, Gini index, mean residual life function (MRLF) also does not exists. The beauty of CEL is that, this is a single parameter decreasing hazard distribution and explains the phenomenon better than others. Some of its statistical properties such as hazard rate function, characteristics function, order statistics, stochastic ordering, quantile function are derived. The expressions of Rényi and Tsallis entropies has also been discussed. The maximum likelihood method of estimation is employed to estimate the parameter and a simulation study is performed to check the consistency of the estimate of the parameter. Although the moments do not exist for CEL, but Figure~\ref{fig:pdf} indicates that, the distribution is highly positively skewed distribution. As the value of $\theta$  is increasing the density of the distribution becomes flatten. Lastly we have considered two different data sets and five other distributions having decreasing failure rate namely exponential Poisson Lindley (EPL), Exponential Poisson (EP),
Exponential Logarithmic (EL), Weibull, gamma distributions. It is
shown that our proposed NWEL distribution fits the considered real data set very well even better than other distributions according to K-S Statistics and associated $p$-value. Hence, we can easily conclude that our proposed CEL distribution may be considered as a suitable model for the case of decreasing failure rate scenario with a hope to get better model in various disciplines such as medical, engineering, and social sciences.

%%----------		References		----------%%


\begin{thebibliography}{00}
%% Journal references
\bibitem[Adamidis \& Loukas(1998)]{adamidis1998lifetime}
        Adamidis, K. \& Loukas, S. A. lifetime distribution with decreasing failure rate.
        Statistics \& Probability Letters, 1998; 39(1): 35-42.       
        
\bibitem[Barlow \& Marshall(1964)]{barlow1964bounds}        
        Barlow, R. E. \& Marshall, A. W. 
        Bounds for distributions with monotone hazard rate. The Annals of Mathematical Statistics, 1964; 35(3), 1234-1257.
        
\bibitem[Barlow \& Marshall(1965)]{barlow1965tables}        
        Barlow, R. E., \& Marshall, A. W.
        Tables of bounds for distributions with monotone hazard rate. 
        Journal of the American Statistical Association, 1965; 60(311), 872-890.
        
\bibitem[Barlow et al.(1963)]{barlow1963properties}        
        Barlow, R. E., Marshall, A. W. \& Proschan, F.
        Properties of probability distributions with monotone hazard rate. The Annals of Mathematical Statistics, 1963; 34(2), 375-389.

\bibitem[Barreto-Souza et al.(2011)]{barreto2011weibull}        
        Barreto-Souza, W., de Morais, A. L. \& Cordeiro, G. M.
        The Weibull-geometric distribution.
        Journal of Statistical Computation and Simulation, 2011; 81(5): 645-657.

\bibitem[Barreto-Souza \& Bakouch(2013)]{barreto2013new}        
        Barreto-Souza, W. \& Bakouch, H. S. 
        A new lifetime model with decreasing failure rate. 
        Statistics, 2013; 47(2): 465-476.
        
\bibitem[Cozzolino(1968)]{cozzolino1968probabilistic}        
        Cozzolino, J. M.
        Probabilistic models of decreasing failure rate processes.
        Naval Research Logistics Quarterly, 1968; 15(3), 361-374.
        
\bibitem[Chahkandi \& Ganjali(2009)]{chahkandi2009some}        
        Chahkandi, M. \& Ganjali, M.
        On some lifetime distributions with decreasing failure rate. Computational Statistics \& Data Analysis, 2009; 53(12): 4433-4440.   
        
\bibitem[Dahiya \& Gurland(1972)]{dahiya1972goodness}        
        Dahiya, R. C. \& Gurland, J.  
        Goodness of fit tests for the gamma and exponential distributions. Technometrics, 1972; 14(3), 791-801. 
        
\bibitem[Ghitany et al.(2011)]{ghitany2011two}
       Ghitany, M. E., Alqallaf, F., Al-Mutairi, D. K. \&  Husain,   H. A.
       A two-parameter weighted Lindley distribution and its applications to survival data. Mathematics and Computers in Simulation, 2011; 81, 1190-1201.

\bibitem[Glaser(1980)]{laser1980bathtub}        
        Glaser, R. E.
        Bathtub and related failure rate characterizations.
        Journal of the American Statistical Association, 1980; 75(371):
        667-672.
        
\bibitem[Gleser(1989)]{gleser1989gamma}        
        Gleser, L. J. The gamma distribution as a mixture of exponential distributions. The American Statistician, 1989; 43(2), 115-117.    
        
\bibitem[Gurland \& Sethuraman(1994)]{gurland1994shorter}        
        Gurland, J. \& Sethuraman, J.
        Reversal of increasing failure rates when pooling failure data. 
        Technometrics, 1994; 36(4), 416-418.
        
\bibitem[Kuş(2007)]{kucs2007new}        
        Kuş, C.
        A new lifetime distribution. 
        Computational Statistics \& Data Analysis, 2007; 51(9): 4497-4509.
        
\bibitem[Lindley(1958)]{lindley1958}
       Lindley D. V.
       Fiducial distributions and Bayes theorem. Journal of the Royal Statistical Society, Series B (Methodological), 1958; 102-107.
        
\bibitem[Linhart \& Zucchini(1986)]{linhart1986}        
        Linhart, H. \& Zucchini, W.
        Model selection. John Wiley \& Sons, 1986. 
        
\bibitem[Lomax(1954)]{lomax1954business}
        Lomax, K. S. Business failures: Another example of the analysis of failure data.
        Journal of the American Statistical Association, 1954; 49(268), 847-852.   
        
\bibitem[Marshall \& Proschan(1965)]{marshall1965maximumy}        
        Marshall, A. W. \& Proschan, F.
        Maximum likelihood estimation for distributions with monotone failure rate. 
        The annals of mathematical statistics, 1965; 36(1), 69-77.            
        
\bibitem[Moors(1988)]{moors1988quantile}
        Moors, J. J. A. 
        A quantile alternative for kurtosis. 
        Journal of the Royal Statistical Society: Series D (The Statistician), 1988; 37(1), 25-32.
        
\bibitem[Nassar(1988)]{nassar1988two}
        Nassar, M. M. 
        Two properties of mixtures of exponential distributions. IEEE transactions on reliability, 1988; 37(4), 383-385.
        
\bibitem[Nelson(1982)]{nelson1982weibull}
        Nelson, W. Applied Life Data Analysis. John Wiley \& Sons, New York, 1982.              
        
\bibitem[Proschan(1963)]{proschan1963theoretical}        
        Proschan, F.  
        Theoretical explanation of observed decreasing failure rate. Technometrics, 1963; 5, 375-383.
        
\bibitem[Rényi(1961)]{renyi1961}        
        Rényi, A.
        On measures of entropy and information. In: Proceedings of the Fourth Berkeley Symposium on Mathematical Statistics and Probability.
        Contributions to the Theory of Statistics, Berkeley, California: University of California Press, 1961; (1): 547–561. 
        
\bibitem[Saunders \& Myhre(1983)]{saunders1983maximum}
        Saunders, S. C. \& Myhre, J. M. 
        Maximum likelihood estimation for two-parameter decreasing hazard rate distributions using censored data. Journal of the American Statistical Association, 1983; 78(383), 664-673.        
        
\bibitem[Shaked \& Shanthikumar(1994)]{shakedjg}        
        Shaked, M. \& Shanthikumar, J.
        Stochastic Orders and Their Applications. Boston: Academic Press, 1994.               
        
\bibitem[Tahmasbi \& Rezaei(2008)]{tahmasbi2008two}        
        Tahmasbi, R. \&  Rezaei, S.
        A two-parameter lifetime distribution with decreasing failure rate. Computational Statistics \& Data Analysis, 2008; 52(8): 3889-3901.
        
\bibitem[Tsallis(1988)]{tsallis1988possible}        
        Tsallis, C. 
        Possible generalization of Boltzmann-Gibbs statistics. Journal of statistical physics, 1988; 52(1-2), 479-487.
                
\bibitem[Weibull(1951)]{weibull1951}        
        Weibull, W. 
        A statistical distribution function of wide applicability. 
        Journal of applied mechanics, 1951; 18(3): 293-297.
        
\end{thebibliography}
\end{document}